\documentclass[a4paper,oneside,12pt]{article}
\usepackage{amsmath,amsfonts,amssymb,amsthm,mathrsfs}
\usepackage[a4paper,vmargin={3.5cm,3.5cm},hmargin={2.5cm,2.5cm}]{geometry}
\usepackage[font=sf, labelfont={sf,bf}, margin=1cm]{caption}
\usepackage{graphicx,graphics}
\usepackage{epsfig}
\usepackage{latexsym}
\usepackage[applemac]{inputenc}
\usepackage{ae,aecompl}
\usepackage[english]{babel}
 \usepackage[colorlinks=true]{hyperref}
\usepackage{pstricks}

\newtheorem{theorem}{Theorem}

\theoremstyle{definition}

\theoremstyle{remark}

\newcommand{\R}{\mathbb R}

\newcommand{\Z}{\mathbb Z}

\renewcommand{\geq}{\geqslant}
\renewcommand{\leq}{\leqslant}

\renewcommand{\leq}{\leqslant}
\renewcommand{\geq}{\geqslant}
\title{\textsc{The Liouville and the intersection properties are equivalent for planar graphs}}
\date{}
\author{Itai Benjamini \and  Nicolas Curien \and Agelos Georgakopoulos}
\begin{document}
\maketitle

\begin{abstract}
It is shown that if a planar graph admits no non-constant bounded harmonic function then the trajectories of two independent simple random walks  intersect almost surely.
\end{abstract}

\section{Introduction}
Let $G=(V,E)$ be a connected (multi)graph. 
The graph $G$ has the \emph{intersection} property if for any $x,y \in V$, the trajectories of two independent simple random walks (SRW) started respectively from $x$ and $y$ intersect almost surely.
Recall that $G$ is {\it Liouville} if it admits no non-constant bounded harmonic function. The goal of this note is to prove the following:
\begin{theorem}\label{main}	
If $G$ is planar and Liouville then $G$ has the intersection property. 
\end{theorem}
Let us make a few comments on this result. We first recall the well-known sequence of implications as well as the $\mathbb{Z}^d, d \geq 1$ lattices that satisfy them:
$$ \begin{array}{ccccc}\mbox{recurrence} & \Rightarrow & \mbox{intersection} & \Rightarrow & \mbox{Liouville}\\
\mathbb{Z}^d, d=1, 2 & & \mathbb{Z}^d, d =1,2,3,4 & & \mathbb{Z}^d, d \geq 1. \end{array}$$

In the case of \emph{bounded degree} planar graphs the Liouville property is equivalent to recurrence of the simple random walk, see \cite{BS96a} (as in the case of non-compact planar Riemannian surfaces). If we drop the bounded degree assumption, it is easy to construct transient planar graphs which are Liouville: For example, start with a half line $  \mathbb{N} = \{ 0,1,2,... \}$ and place $2^n$ parallel edges between $n$ and $n+1$ for $n \geq 0$. Yet this graph clearly has the intersection property.


Our result thus becomes interesting in the case of Liouville planar graphs of unbounded degree. Such graphs arise for example as distributional limits of finite random planar graphs, a topic that has attracted a lot of research recently \cite{AS03,BCstationary,GGN12,Kri05}.


\section{Proof}
\begin{figure}[!h]
  \begin{center}
 \includegraphics[width=8cm]{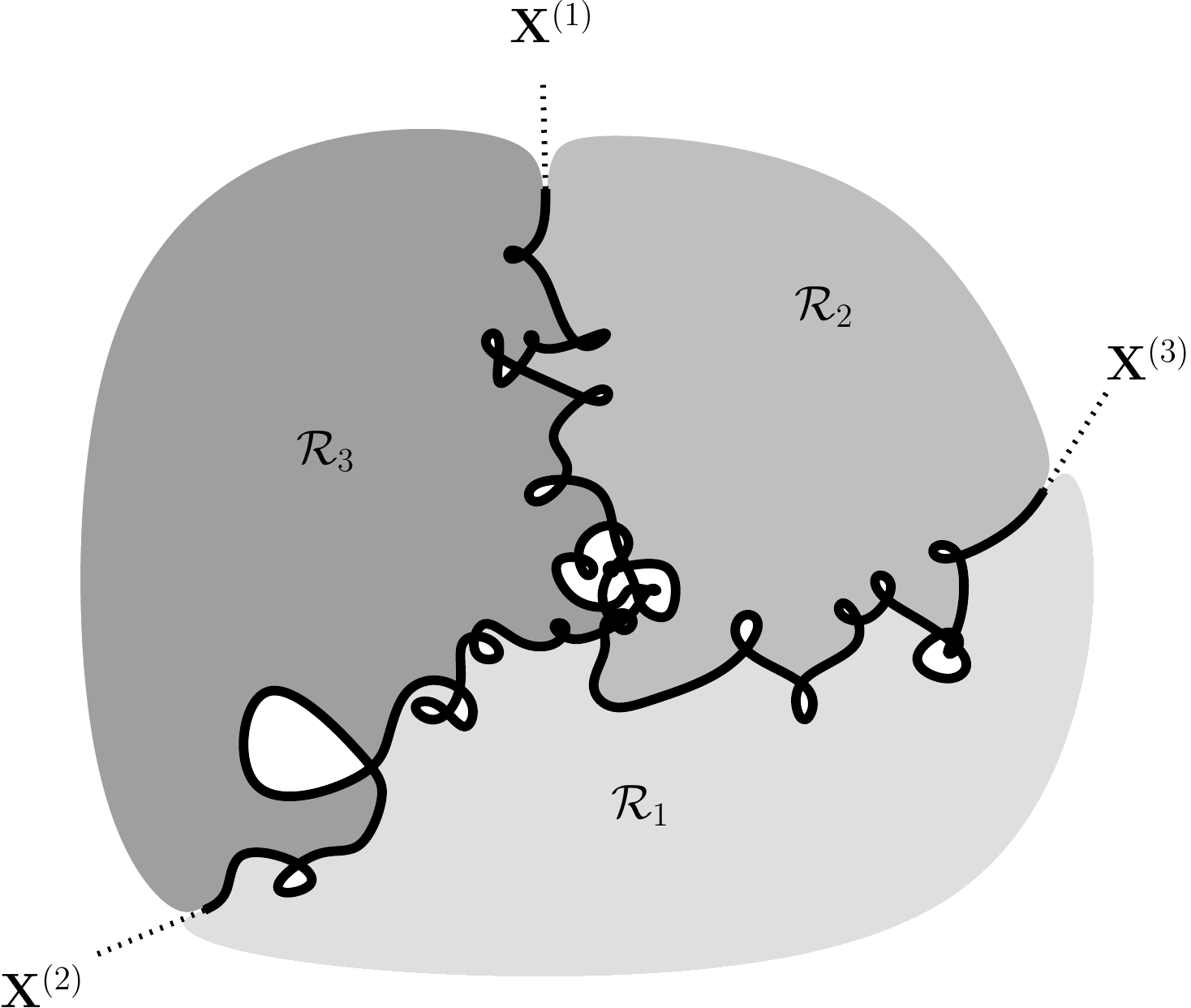}
  \caption{ \label{separation}The three (random) regions $ \mathcal{R}_{1}, \mathcal{R}_{2}$ and $ \mathcal{R}_{3}$.}
  \end{center}
  \end{figure}

The strategy of the proof is the following. We consider three independent simple random walk trajectories, and argue that if each two of them intersect only finitely many times, then they divide our planar graph into three regions (Fig.\ref{separation}). This allows us to talk about the probability for random walk to eventually stay in one of these regions, which we use to construct a non-constant bounded harmonic function, contradicting our Liouvilleness assumption.

\proof Fix $G=(V,E)$ a connected planar multi-graph and suppose that $G$ has the Liouville property. First note that if $G$ is recurrent then it has the intersection property and we thus suppose henceforth that $G$ is transient. \\
Denote by $P_x$ the law of a simple random walk  $ (X_{n})_{n \geq 0}$ started from $x$ in $G$. We write $\mathbf{X} = \{ X_{0}, X_{1}, X_{2}, ... \}$ for the trajectory of $(X_{n})$. If $\gamma = \{\gamma_0,\gamma_1,...\}$ is a set of vertices in $G$ define for any $x \in V$
 \begin{eqnarray*}h_\gamma(x) &:=& P_x\big( \#(\mathbf{X} \cap \gamma) = \infty\big).  \end{eqnarray*}
 It is clear that $h_\gamma(.)$ is harmonic (bounded) and is thus constant since $G$ is Liouville. We write $h_{\gamma}$ for the common value $h_{\gamma}(y),$ $y \in V$. In fact, one has $h_\gamma \in \{0,1\}$. Indeed, if $ \mathcal{F}_n = \sigma(X_0,X_1, ... , X_n)$ is the sigma-field generated by the SRW we have  \begin{eqnarray*}
\mathbf{1}_{ \# (\mathbf{X} \cap \gamma) =\infty} &=& \lim_{n \to \infty} E_x[\mathbf{1}_{ \# (\mathbf{X} \cap \gamma) =\infty} \mid \mathcal{F}_n] \quad = \quad \lim_{n\to \infty}h_\gamma(X_n)\qquad a.s.  \end{eqnarray*} Hence $ h_{\gamma}(X_{n})$ tends to $0$ or $1$ a.s. and the claim follows.  Let us now randomize the path $\gamma$ and consider the random variable $h_ \mathbf{X}$ where $ \mathbf{X}$ is the trajectory of a SRW under $P_{x}$. The last argument shows that $h_ \mathbf{X} \in \{0,1\}$. We now claim that $h_ \mathbf{X}$ is almost surely constant. Indeed, if $ \mathbf{X}$ and $ \mathbf{Y}$ are two independent simple random walk trajectories started from $x \in V(G)$  we have $ h_{ \mathbf{X}} = h_{ \mathbf{X}}(x) = \mathbf{1}_{\#( \mathbf{X} \cap \mathbf{Y})= \infty}$ a.s. (and similarly interchanging the roles of $ \mathbf{X}$ and $ \mathbf{Y}$) thus
 \begin{eqnarray*} E_x[h_ \mathbf{X}h_ \mathbf{Y}] &=& P_x\big( \# (\mathbf{X} \cap \mathbf{Y}) = \infty \mbox{ and } \# (\mathbf{Y} \cap \mathbf{X})= \infty\big)\\  &=& P_x\big( \# (\mathbf{X} \cap \mathbf{Y}) = \infty\big)\\
 &=& E_x[h_ \mathbf{X}].\\  \end{eqnarray*} Hence $E_x[h_ \mathbf{X}]^2 = E_x[h_ \mathbf{X}] \in \{0,1\}$. Either $h_ \mathbf{X}=1= h_{ \mathbf{X}}(y)$ for all $y \in V$ a.s. in which case Theorem \ref{main} is proved, or almost surely for all $y \in V$ we have $h_ \mathbf{X}=0= h_{ \mathbf{X}}(y)$. In words, a.s. for any $x,y \in V$, two simple random walk paths started from $x$ and $y$ intersect only finitely many times. Let us suppose by contradiction that we are in the latter case.\\

 
 We now make use of the planarity and consider a proper embedding\footnote{Notice that since $G$ is Liouville, it has precisely one transient end $\alpha$ and it is possible to embed $G$ in $ \mathbb{R}^2$ so that no ray in $\alpha$ has an accumulation point in $\R^2$, which always exists \cite{RT02}.} of $G$ in $\R^2$.  Let us  consider three independent simple random walk trajectories $ \mathbf{X}^{(1)}, \mathbf{X}^{(2)}$ and $\mathbf{X}^{(3)}$ started from some $x \in V$. Since almost surely these paths intersect each-other  finitely often a.s., they distinguish, using the planarity, three regions in (the embedding of) $G$ called $ \mathcal{R}_{1}, \mathcal{R}_{2}$ and $ \mathcal{R}_{3}$ as depicted in the figure below.  Formally, the region $ \mathcal{R}_{1}$ is defined as the set of all vertices $y \in  V \backslash \bigcup_{k =1}^3 \mathbf{X}^{(k)}$ such that for all  $n$ large enough, if $\gamma$ is a path linking $y$ to $X_{n}^{(1)}$ in $G$ then $\gamma$ must intersect $ \mathbf{X}^{(2)} \cup \mathbf{X}^{(3)}$. The regions  $ \mathcal{R}_{2}$ and $ \mathcal{R}_{3}$ are defined similarly.
 
  Now for any $y \in V$, conditionally on $ \mathbf{X}^{(1)}, \mathbf{X}^{(2)}$ and $\mathbf{X}^{(3)}$ we consider a simple random walk trajectory $ \mathbf{Y}$ under $P_{y}$. By our assumption the path $ \mathbf{Y}$ intersects any of the $ \mathbf{X}^{(i)}$'s finitely many times a.s. so that $ \mathbf{Y}$ is eventually trapped in one of the three regions $ \mathcal{R}_{1}, \mathcal{R}_{2}$ or $ \mathcal{R}_{3}$. We then define  \begin{eqnarray*} \mathcal{H}(y) &:= & P_{y}( \mathbf{Y}  \mbox{ is eventually trapped in } \mathcal{R}_1).  \end{eqnarray*}
  Note that $ \mathcal{H}(.)$ is a random function and that for every $\omega$ it is (bounded) harmonic over $G$ and thus constant by the Liouville property of $G$.  On the one hand an obvious symmetry argument shows that $ E[ \mathcal{H}(x)] = 1/3$. On the other hand since $ \mathcal{H}$ is almost surely constant we have $ \mathcal{H}(x) = \mathcal{H}(X_{n}^{(1)})$ for all $n \geq 0$. Almost surely, for all  $n$ large enough, if a simple random walk started from $X_{n}^{(1)}$ is eventually trapped in $ \mathcal{R}_{1}$ then it has to cross one of the paths $ \mathbf{X}^{(2)}$ or $ \mathbf{X}^{(3)}$. We thus have \begin{eqnarray*} E[ \mathcal{H}(x)] & = & \limsup_{n \to \infty}  E[\mathcal{H}(X_{n}^{(1)})] \\ &= & \limsup_{n \to \infty} E\Big[P_{X_{n}^{(1)}}( \mathbf{Y} \mbox{ is trapped in } \mathcal{R}_{1})\Big] \\ & \leq & \limsup_{n\to \infty } E \Big[ P_{X_{n}^{(1)}}( \mathbf{Y} \mbox	{ intersects } \mathbf{X}^{(2)} \cup \mathbf{X}^{(3)}) \Big]. \end{eqnarray*}
But since $ \mathbf{X}^{(1)}$ has only a finite intersection with $ \mathbf{X}^{(2)} \cup \mathbf{X}^{(3)}$ a.s. we deduce that the probability inside the expectation of the right-hand side of the last display goes to $0$ as $n \to \infty$. Hence $ E[\mathcal{H}(x)]=0$ and this is a contradiction. One deduces that $ h_{ \mathbf{X}} =1$ almost surely as desired. \endproof

\section{Concluding Remarks}

\paragraph{Relaxing planarity.}  As we noted above, $\Z^5$ is Liouville and yet two independent simple random walk paths do not intersect with positive probability. For which families of graphs, other then planar, intersection is equivalent to Liouville?  E.g., a natural extension of the collection of all planar graphs is the collection
of all graphs with an excluded minor. Fix a finite graph  $H$. It is reasonable to guess that the statement of Theorem \ref{main} still holds if we replace the planarity assumption by the hypothesis that $G$ does not have $H$ as a minor. Indeed, many theorems on planar graphs generalize to excluded minor graphs (see,
e.g., \cite{Tho99}). It is even more interesting to show that  bounded degree transient
excluded minor graphs admit non-constant bounded harmonic functions thus extending \cite{BS96a}.  Another generalization of planarity is sphere-packability, see \cite{BS09}. We \emph{ask} whether a Liouville sphere-packable graph in $ \mathbb{R}^3$ should admit the intersection property?

\paragraph{Quasi-isometry. } Let $G_{1}=(V_{1},E_{1})$ and $G_{2}=(V_{2},E_{2})$ be two graphs endowed respectivelly with their graph distances $d_{1}$ and $d_{2}$. These graphs are called {\em quasi-isometric} if there exists $\phi : V_{1} \to V_{2}$ and two constants $A,B>0$ such that for every $x,y \in V_{1}$ $$A^{-1}d_{1}(x,y)-B \leq d_{2}(\phi(x),\phi(y)) \leq A d_{1}(x,y) +B, $$ and if for every $z \in V_{2}$ there exists $x \in V_{1}$ with $d_{2}(z,\phi(x)) \leq B$.  \\
For bounded degree planar graphs, the Liouville property is quasi-isometry invariant since transience is. An example in \cite{BS96a} shows that Liouville property and almost sure intersection of simple random walks are not quasi-isometry invariants:   there exist two graphs $G_{1}$ and $G_{2}$ that are quasi-isometric such that $G_{1}$ has the intersection property (hence Liouville) whereas $G_{2}$ is non Liouville and two independent simple random walk trajectories started from two different points in $G_{2}$ have a probability strictly in between $0$ and $1$ of having an infinite intersection. However it is possible to modify and iterate the construction of \cite{BS96a} in such a way that the last probability is actually $0$.

\paragraph{Finite Dirichlet energy.} 
A concept related to Liouvilleness is the {\em (Dirichlet) energy} of a harmonic function $f$, defined by
$ \sum_{x \sim y} |f(x)-f(y)|^2$ where the sum ranges over all neighbors $x,y$ of the graph. For example, it is known that if a graph admits a non-constant  harmonic function of finite energy then it is not Liouville \cite{Soardi}, but the converse is in general not true. Russ Lyons asked (private communication) whether the converse becomes true for planar graphs, i.e.\ whether there is a planar graph which admits non-constant bounded harmonic functions but yet none of finite energy. We construct such a graph here. Start with the integers $  \mathbb{Z} = \{ 0,1,2,... \}$, and place $2^{|n|}$ parallel edges between $n$ and $n+1$ for every $n$. Then, join each $n$ to $-n$ by a new edge of resistance $n$, or equivalently, with a path of length $n$ of unit resistance edges. This graph is, easily, still non-Liouville, and the reader will be able to check that it does however admit no non-constant  harmonic functions of finite energy, see for example \cite[Lemma 3.1.]{lamps}.


\end{document}